\documentclass[11.5pt]{article}
\usepackage{amsmath,amssymb,color,fancyhdr}
 \allowdisplaybreaks  

\def\R{{\mathbb{R}}}

\def\BMO{{\mathrm{BMO}}}
\def\loc{{\mathrm{loc}}}

\usepackage{graphicx,amssymb,mathrsfs,amsmath,latexsym,amsfonts,titlesec,amscd,epsfig,cite,amsthm}
\usepackage{indentfirst} 

\DeclareMathOperator*{\essinf}{ess\, inf}
\DeclareMathOperator*{\esssup}{ess\, sup}
\DeclareMathOperator*{\supp}{supp}

\newcommand{\upcite}[1]{\textsuperscript{\textsuperscript{\cite{#1}}}}

\theoremstyle{definition} 

\newtheorem{theorem}{\indent  Theorem}[section]
    \newtheorem{lemma}{\indent  Lemma} [section]

\newtheorem{proposition}[theorem]{\indent  Proposition}

\theoremstyle{definition} 
    \newtheorem{definition}{\indent  Definition} [section]

\newtheorem*{acknowledgments}{\indent Acknowledgments}
\title{\bf\Large
 Boundedness for commutators of fractional integrals  on  Herz-Morrey spaces with variable exponent
  }

\author{
{\normalsize Jianglong WU
}
\\
{\small\it   Department of Mathematics,  Mudanjiang Normal University, Mudanjiang \ 157011, China}\\
  }
\date{} 

\begin{document}

\maketitle



\begin{minipage}[t]{14cm}

\setlength{\baselineskip}{1.0em}

\noindent

 { \bf Abstract:} In this paper, some boundedness for commutators of fractional integrals are obtained on  Herz-Morrey spaces with variable exponent applying some properties of varible exponent and $\BMO$ function.


 \smallskip
 { \bf AMS(2010) Subject Classification:}  \ 47B47; 42B20; 42B35   \ \

\end{minipage}

\section{Introduction}

Function spaces with variable exponent are being watched with keen interest not in real analysis but also in partial differential equations and in applied mathematics because they are applicable to the modeling for electrorheological fluids and image restoration. The theory of function spaces with variable exponent has rapidly made progress in the past twenty years
since some elementary properties were established by Kov\'{a}\v{c}ik and R\'{a}kosn\'{i}k\upcite{KR}. One of the main
problems on the theory is the boundedness of the Hardy-Littlewood maximal operator on
variable Lebesgue spaces. By virtue of the fine works\upcite{CDF,CFMP,CFN,D1,D2,DHHMS,K,L,N,PR}, some
important conditions on variable exponent, for example, the $\log$-H\"{o}lder conditions et al, have been obtained.

The class of the Herz spaces is arising from the study on characterization of multipliers on the classical Hardy spaces. And the homogeneous Herz-Morrey spaces $M\dot{K}^{\alpha,\lambda}_{p,q}(\R^{n})$ coordinate with  the homogeneous Herz space $\dot{K}^{\alpha,p}_{q}(\R^{n})$ when $\lambda=0$. One of the important problems on Herz spaces and Herz-Morrey spaces is the boundedness of sublinear operators. Hern\'{a}ndez, Li, Lu and Yang et al\upcite{HY,LiY,LY} have proved that if a sublinear operator $T$ is bounded on $L^{p}(\R^{n})$ and satisfies the size condition
$$ |Tf(x)|\le C \int_{\R^{n}}\frac{|f(y)|}{|x-y|^{n}}\mathrm{d}y$$
for all $f\in L^{1}(\R^{n})$ with compact support and a.e. $x\notin\supp f$, then $T$ is bounded on the homogeneous Herz space $\dot{K}^{\alpha,p}_{q}(\R^{n})$. In 2005, Lu and Xu\upcite{LX} established the boundedness for some sublinear operators.

The $\BMO$ space and the $\BMO$ norm are defined respectively as follows:
\begin{equation*}
\begin{split}
\BMO(\R^{n}) = \Big\{ b\in L^{1}_{\loc}(\R^{n}) :\|b\|_{\BMO(\R^{n})}<\infty \Big\},    \ \
\|b\|_{\BMO(\R^{n})} = \sup_{B: ball} \frac{1}{|B|} \int_{B}|b(x)-b_{B}| \mathrm{d}x.
\end{split}
\end{equation*}
The fractional integral $I_{\beta}$ is defined by
$I_{\beta}(f)(x) = \int_{\R^{n}}\frac{f(y)}{|x-y|^{n-\beta}} \mathrm{d}y,
$
the commutator for fractional integral is defined by
$[b,I_{\beta}]f(x)=b(x)I_{\beta}(f)(x)-I_{\beta}(bf)(x),$
and $m$-order commutator for fractional integral is defined by
$$I_{\beta,b}^{m}(f)(x)=\int_{\R^{n}}\frac{f(y)(b(x)-b(y))^{m}}{|x-y|^{n-\beta}} \mathrm{d}y,$$
where $0<\beta<n, b\in \BMO(\R^{n}), m\in \mathbb{N}$. It is easy to see, when $m=1$, $I_{\beta,b}^{m}(f)(x)=[b,I_{\beta}]f(x)$; and when $m=0$, $I_{\beta,b}^{m}(f)(x)=I_{\beta}(f)(x)$.

Chanillo\upcite{C} has initially introduced the commutator $[b,I_{\beta}]$ with $b\in \BMO$ and proved the boundedness on Lebesgue spaces with constant exponent. In 2010, Izuki\upcite{I3} generalizes this result to the case of variable exponent and  considere the boundedness on Herz spaces with variable exponent.

In 2010, Izuki\upcite{I2} proves the boundedness of some sublinear operators on Herz spaces with variable exponent. And recently Izuki\upcite{I,I1} also considers the boundedness of some operators on Herz-Morrey  spaces with variable exponent.

Motivated by the study on the Herz spaces and Lebesgue spaces  with variable exponent, the main purpose of this paper is to establish some boundedness for commutators of fractional integrals on Herz-Morrey spaces with variable exponent, Our main tools are some properties of varible exponent and $\BMO$ function. And we also note that our results are the generalizations of main theorems for Izuki\upcite{I,I3} on Herz space and Herz-Morrey spaces with variable exponent.

Throughout this paper, we will denote by $|S|$ the Lebesgue measure  and  by $\chi_{_{\scriptstyle S}}$ the characteristic function  for a measurable set  $S\subset\R^{n}$. Given a function $f$, we denote the mean value of $f$ on $S$ by $f_{S}:=\frac{1}{|S|} \int_{S}f(x)\mathrm{d}x$.  $C$  denotes a constant that is independent of the main parameters involved but
whose value may differ from line to line. For any index $1< q(x)< \infty$, we denote by $q'(x)$ its conjugate index,
namely, $q'(x)=\frac{q(x)}{q(x)-1}$.  For $A\sim D$, we mean that there is a constant $C > 0$ such that$C^{-1}D\le A \le CD$.

\section{Preliminaries  and Lemmas}

In this section, we give the definition of Lebesgue and Herz-Morrey spaces with variable exponent, and
state their properties. Let $E$ be a measurable set in $\R^{n}$ with $|E|>0$. We first define Lebesgue spaces with variable exponent.

\begin{definition} \label{Def.L}
Let ~$ q(\cdot): E \to[1,\infty)$ be a measurable function.

1)\ The Lebesgue spaces with variable exponent $L^{q(\cdot)}(E)$ is defined by
  $$ L^{q(\cdot)}(E)=\{f~ \mbox{is measurable function}: \int_{E} \Big( \frac{|f(x)|}{\eta} \Big)^{q(x)} \mathrm{d}x <\infty ~\mbox{for some constant}~ \eta>0\}.  $$

  2)\ The space $L_{\loc}^{q(\cdot)}(E)$ is defined by
  $$ L_{\loc}^{q(\cdot)}(E)=\{f ~\mbox{is measurable function}: f\in L^{q(\cdot)}(K) ~\mbox{for all compact  subsets}~ K\subset E\}.  $$
\end{definition}

The Lebesgue space $L^{q(\cdot)}(E)$ is a Banach space with the norm defined by  $$ \|f\|_{L^{q(\cdot)}(E)}=\inf \Big\{ \eta>0:  \int_{E} \Big( \frac{|f(x)|}{\eta} \Big)^{q(x)} \mathrm{d}x \le 1 \Big\}.$$

Now, we define two classes of exponent functions. Given a function $f\in L_{\loc}^{1}(E)$, the Hardy-Littlewood maximal operator $M$ is defined by
$$Mf(x)=\sup_{r>0} r^{-n} \int_{B(x,r)\cap E} |f(y)| \mathrm{d}y \ \ \ (x\in E),$$
where $B(x,r)=\{y\in \R^{n}: |x-y|<r\}$.

\begin{definition} \label{Def.PB}
1)\ The set $\mathscr{P}(\R^{n})$ consists of all measurable functions $ q(\cdot)$ satisfying
$$1< \essinf_{x\in \R^{n}} q(x)=q_{-},\ \ q_{+}= \esssup_{x\in \R^{n}} q(x)<\infty.$$

2)\  The set $\mathscr{B}(\R^{n})$ consists of all  measurable functions  $q(\cdot)\in\mathscr{P}(\R^{n})$ satisfying that the Hardy-Littlewood maximal operator $M$ is bounded on $L^{q(\cdot)}(\R^{n})$.
\end{definition}

Next we define the Herz-Morrey spaces with  variable exponent. Let $B_{k}=B(0,2^{k})=\{x\in\R^{n}:|x|\leq 2^{k}\}, A_{k}=\ B_{k}\setminus  B_{k-1}$ and $\chi_{_{k}}=\chi_{_{A_{k}}}$ for $k\in \mathbb{Z}$.

\begin{definition} \label{Def.M-H}
 Let $ \alpha\in \mathbb{R}, ~0\leq \lambda < \infty,~ 0<p< \infty$,
and $q(\cdot)\in \mathscr{P}(\mathbb{R}^{n})$. The Herz-Morrey space with  variable exponent $ M\dot{K}_{p, q(\cdot)}^{\alpha,
\lambda}(\mathbb{R}^{n} ) $ is definded by
  $$ M\dot{K}_{p, q(\cdot)}^{\alpha, \lambda}(\mathbb{R}^{n})=\Big\{f\in L_{\loc}^{q(\cdot)}(\mathbb{R}^{n}\backslash\{0\}):
 \|f\|_{M\dot{K}_{p, q(\cdot)}^{\alpha, \lambda}(\mathbb{R}^{n})}<\infty \Big\},  $$
where
  $$ \|f\|_{M\dot{K}_{p, q(\cdot)}^{\alpha, \lambda}(\mathbb{R}^{n})}=\sup_{k_{0}\in \mathbb{Z}}2^{-k_{0}\lambda}
 \Big(\sum_{k=-\infty}^{k_{0}}2^{k\alpha
 p}\|f\chi_{_{\scriptstyle k}}\|_{L^{^{q(\cdot)}}(\mathbb{R}^{n})}^{p} \Big)^{\frac{1}{p}}.  $$

\end{definition}

Compare the Herz-Morrey space with  variable exponent  $ M\dot{K}_{p, q(\cdot)}^{\alpha, \lambda}(\mathbb{R}^{n})$ with the Herz space with  variable exponent $\dot{K}_{q(\cdot)}^{\alpha,p}(\R^{n})$\upcite{I1}, where
$$\dot{K}_{q(\cdot)}^{\alpha,p}(\R^{n})= \Big\{f\in L_{\loc}^{q(\cdot)}(\mathbb{R}^{n}\backslash\{0\}):
\sum\limits_{k=-\infty}^{\infty}2^{k\alpha  p}\|f\chi_{_{k}}\|_{L^{q(\cdot)}(\R^{n})}^{p}<\infty \Big\},$$
Obviously, $M\dot{K}_{p,q(\cdot)}^{\alpha,0} (\R^{n})=\dot{K}_{q(\cdot)}^{\alpha,p}(\R^{n})$.

 When $\lambda=0$, we can see that our result below generalize the result in the setting of the Herz space with variable exponent, which proved by  Izuki in \cite{I3}. So in this paper, we only give the result when $\lambda>0$.

In 2012,  Almeida and  Drihem\upcite{AD} discuss the boundedness of a wide class of sublinear operators, including maximal, potential and Calder\'{o}n-Zygmund operators, on variable Herz spaces $K_{q(\cdot)}^{\alpha(\cdot),p}(\R^{n})$ and $\dot{K}_{q(\cdot)}^{\alpha(\cdot),p}(\R^{n})$. Meanwhile, they also establish Hardy-Littlewood-Sobolev theorems for fractional integrals on variable Herz spaces. In this paper, the author only considers Herz-Morrey space $ M\dot{K}_{p, q(\cdot)}^{\alpha(\cdot), \lambda}(\mathbb{R}^{n})$ with variable exponent $q(\cdot)$ but fixed $ \alpha\in \mathbb{R}$ and $p\in(0,\infty)$. However, for  the case of the exponent $\alpha(\cdot)$ is variable as well, we can refer to the furthermore work for the author.

Next we state some properties of variable exponent. Cruz-Uribe et al\upcite{CFN} and Nekvinda\upcite{N} proved the following sufficient conditions independently. Moreover, we note that Diening\upcite{D2} proved the following proposition in the case of $E$ is bounded, and Nekvinda\upcite{N} gave a more general condition in place of (\ref{C2}).

\begin{proposition}\ \label{Pr-1}
Suppose that $E$ is an open set, If $ q(\cdot)\in \mathscr{P}(E)$ satisfies the inequality
 \begin{alignat}{2}
 |q(x)-q(y)| &\le \frac{-C}{\ln(|x-y|)}  \ \ \ \mbox{if}~ |x-y|\le 1/2,   \label{C1}\\
 |q(x)-q(y)| &\le \frac{C}{\ln(e+|x|)}  \ \ \ \mbox{if}~ |y|\ge|x|,\ \    \label{C2}
\end{alignat}
where $C>0$ is a constant independent of $x$ and $y$, then we have $ q(\cdot)\in \mathscr{B}(E)$.
\end{proposition}

In order to prove our main theorem, we also need the following result which is the Hardy-Littlewood-Sobolev theorem on Lebesgue spaces with varible expoonent due to Capone, Cruz-Uribe and Fiorenza\upcite{CCF}(see Theorem 1.8). We remark that this result is initially proved by Diening\upcite{D} provided that $q_{_{1}}(\cdot)$ is constant outside of a large ball.

\begin{proposition}\upcite{CCF}\ \label{Pr-3}
Suppose that $q_{_{1}}(\cdot)\in \mathscr{P}(\R^{n})$ satisfies conditions (\ref{C1}) and (\ref{C2}) in Proposition \ref{Pr-1}.
$0<\beta<n/(q_{_{1}})_{+}$ and define  $q_{_{2}}(\cdot)$ by
$$\frac{1}{q_{_{1}}(x)} -\frac{1}{q_{_{2}}(x)}=\frac{\beta}{n}.
$$
Then we have
$$\|I_{\beta}f\|_{L^{q_{2}(\cdot)}(\R^{n})} \le {C}\|f\|_{L^{q_{1}(\cdot)}(\R^{n})}$$
for  all $f\in{L^{q_{_{1}}(\cdot)}(\R^{n})}$.
\end{proposition}

In addition, the following result for the boundedness of $I_{\beta,b}^{m}$   on the Lebesgue spaces with variable exponent will be used in the proof of our main theorem.

\begin{proposition}\ \label{Pr-2}
Suppose that $q_{_{1}}(\cdot)\in \mathscr{P}(\R^{n})$ satisfies conditions (\ref{C1}) and (\ref{C2}) in Proposition \ref{Pr-1}. Let $m\in \mathbb{N}, ~ 0<\beta<n/(q_{_{1}})_{+}$, Define the variable exponent $q_{_{2}}(\cdot)$ by
$$\frac{1}{q_{_{1}}(x)} -\frac{1}{q_{_{2}}(x)} = \frac{\beta}{n}.$$
  Then $I_{\beta,b}^{m}$ is  bounded from $L^{q_{_{1}}(\cdot)}(\R^{n})$ into $L^{q_{_{2}}(\cdot)}(\R^{n})$ for all $f\in L^{q_{_{1}}(\cdot)}(\R^{n})$ and $b\in \BMO(\R^{n})$.
\end{proposition}

The idea of the proof for Proposition \ref{Pr-2} comes from the Theorem 1 in \cite{I3}. We omit the details.

The next lemma describes the generalized H\"{o}lder's inequality and the
duality of $L^{q(\cdot)}(E)$. The proof is found in \cite{KR}.

\begin{lemma}\upcite{KR}\  \label{Lem.1}
Suppose that $ q(\cdot)\in \mathscr{P}(E)$, Then the following statements hold.

1)\  (generalized H\"{o}lder's inequality)\ \ For all $f\in L^{q(\cdot)}(E)$ and all $g\in L^{q'(\cdot)}(E)$, we have
 \begin{equation*}
 \int_{E} |f(x)g(x)|\mathrm{d}x \le r_{q}\|f\|_{L^{q(\cdot)}(E)} \|g\|_{L^{q'(\cdot)}(E)},
\end{equation*}
where $r_{q}=1+1/q_{-}-1/q_{+}$.

2)\ For all $f\in L^{q(\cdot)}(E)$, we have
$$\|f\|_{L^{q(\cdot)}(E)} \le \sup \Big\{ \int_{E}|f(x)g(x)|\mathrm{d}x: \|g\|_{L^{q'(\cdot)}(E)}\le 1 \Big\}.$$

\end{lemma}

\begin{lemma}\upcite{I}\ \label{Lem.2}
If $ q(\cdot)\in \mathscr{B}(\R^{n})$, then there exists a positive constant $\delta\in(0,1)$ and $C>0$ such that
 \begin{equation*}
 \frac{\|\chi_{S}\|_{L^{q(\cdot)}(\R^{n})}}{\|\chi_{B}\|_{L^{q(\cdot)}(\R^{n})}} \le C \Big(\frac{|S|}{|B|} \Big)^{\delta}
\end{equation*}
holds for all balls $B$ in $\R^{n}$ and all measurable subsets $S\subset B$.
\end{lemma}

\begin{lemma}\upcite{I}\  \label{Lem.3}
If $ q(\cdot)\in \mathscr{B}(\R^{n})$, then there exists a positive constant  $C>0$ such that
 \begin{equation*}
 C^{-1}\le \frac{1}{|B|} \|\chi_{B}\|_{L^{q(\cdot)}(\R^{n})} \|\chi_{B}\|_{L^{q'(\cdot)}(\R^{n})} \le C
\end{equation*}
for all balls $B$ in $\R^{n}$.
\end{lemma}

\begin{lemma}\upcite{I2}\  \label{Lem.4}
Let $b\in\BMO(\R^{n}),  m\in \mathbb{N},~i,j\in \mathbb{Z}$ with $i<j$. Then we have
  \begin{alignat*}{2}
 C^{-1}\|b\|^{m}_{\BMO(\R^{n})} &\le \sup_{B} \frac{1}{\|\chi_{B}\|_{L^{q(\cdot)}(\R^{n})}}  \|(b-b_{B})^{m}\cdot \chi_{B}\|_{L^{q(\cdot)}(\R^{n})} \le C\|b\|^{m}_{\BMO(\R^{n})},       \\
\|(b-b_{B_{i}})^{m} &\cdot \chi_{B_{j}}\|_{L^{q(\cdot)}(\R^{n})} \le C(j-i)^{m}\|b\|^{m}_{\BMO(\R^{n})} \|\chi_{B_{j}}\|_{L^{q(\cdot)}(\R^{n})}.
\end{alignat*}
\end{lemma}

The above result is proved by Izuki\upcite{I2}. We remark that Lemma \ref{Lem.4} is a generalization of well-known properties for $\BMO$ spaces.

\section{Main theorem and its proof}

In this section we prove the boundedness for the higher order commutator of fractional integrals
on Herz-Morrey spaces with variable exponent under some conditions.

Let $ q(\cdot)\in \mathscr{P}(\R^{n})$ satisfy conditions (\ref{C1}) and (\ref{C2}) in Proposition \ref{Pr-1}. Then so does $ q'(\cdot)$. In particular, we can see that  $ q(\cdot),~q'(\cdot)\in \mathscr{B}(\R^{n})$ from  Proposition \ref{Pr-1}. Therefore applying Lemma \ref{Lem.2} when $q_{1}(\cdot), q_{2}(\cdot)\in \mathscr{P}(\R^{n})$, we can take  constant $\delta_{1}\in (0, 1/(q'_{2})_{+}), \delta_{2}\in (0, 1/(q_{1})_{+})$ such that
\begin{equation}  \label{BZ}
 \frac{\|\chi_{S}\|_{L^{q'_{1}(\cdot)}(\R^{n})}} {\|\chi_{B}\|_{L^{q'_{1}(\cdot)}(\R^{n})}} \le C \bigg(\frac{|S|}{|B|} \bigg)^{\delta_{1}},
 \ \ \  \frac{\|\chi_{S}\|_{L^{q_{2}(\cdot)}(\R^{n})}} {\|\chi_{B}\|_{L^{q_{2}(\cdot)}(\R^{n})}} \le C \bigg(\frac{|S|}{|B|} \bigg)^{\delta_{2}}
\end{equation}
for all balls $B$ in $\R^{n}$ and all measurable subsets $S\subset B$.

 Our main result can be stated as follows.

\begin{theorem}\label{thm.1}
Suppose that $q_{_{1}}(\cdot)\in \mathscr{P}(\R^{n})$ satisfies conditions (\ref{C1}) and (\ref{C2}) in Proposition \ref{Pr-1}. Define the variable exponent $q_{_{2}}(\cdot)$ by
$$\frac{1}{q_{_{1}}(x)} -\frac{1}{q_{_{2}}(x)} = \frac{\beta}{n}.$$
Let $m\in \mathbb{N}, ~ 0<p_{_{1}}\le {p_{_{2}}} <\infty,~ \lambda>0,~ 0<\beta<n/(q_{_{1}})_{+},~ \lambda-n\delta_{2}<\alpha<\lambda+ n\delta_{1}$,  where $\delta_{1}\in (0, 1/(q'_{1})_{+})$ and $\delta_{2}\in (0, 1/(q_{2})_{+})$ are the constants appearing in (\ref{BZ}).
Then  $I_{\beta,b}^{m}$ is  bounded from $M\dot{K}_{p_{_{1}},q_{_{1}}(\cdot)}^{\alpha,\lambda}(\R^{n})$ into $M\dot{K}_{p_{_{2}},q_{_{2}}(\cdot)}^{\alpha,\lambda}(\R^{n})$ for all $f\in M\dot{K}_{p_{_{1}},q_{_{1}}(\cdot)}^{\alpha,\lambda}(\R^{n})$ and $b\in \BMO(\R^{n})$.
\end{theorem}

     \begin{proof}
For $\forall~ f\in {M\dot{K}_{p_{_{1}},  q_{_{1}}(\cdot)}^{\alpha,\lambda}(\R^{n})}$ and $\forall~b\in \BMO(\R^{n})$. If we denote   $f_j:=f\cdot\chi_{j}=f\cdot\chi_{A_j}$ for each $j\in \mathbb{Z}$, then we can write
$$f(x)=\sum_{j=-\infty}^{\infty}f(x)\cdot\chi_{j}(x) =\sum_{j=-\infty}^{\infty}f_{j}(x).
$$
Because of $0<p_{_{1}}/p_{_{2}}\le 1$, we apply inequality
\begin{equation*}
\bigg(\sum_{i=-\infty}^{\infty}|a_{i}|\bigg)^{ p_{_{1}}/p_{_{2}}} \le \sum_{i=-\infty}^{\infty} |a_{i}|^{ p_{_{1}}/ p_{_{2}}},
\end{equation*}
and obtain
\begin{eqnarray*}
&&\; \|I_{\beta, b}^{m}(f)\|^{p_{_{1}}}_{M\dot{K}_{p_{_{2}}, q_{_{2}}(\cdot)}^{\alpha,\lambda}(\R^{n})}= \sup_{k_{0}\in \mathbb{Z}}2^{-k_{0}\lambda {p_{_{1}}}}  \bigg(\sum_{k=-\infty}^{k_{0}}2^{k\alpha {p_{_{2}}}}  \| I_{\beta, b}^{m}(f)\cdot\chi_{_{\scriptstyle k}}  \|_{L^{q_{_{2}}(\cdot)}(\R^{n})}^{p_{_{2}}}\bigg)^{p_{_{1}}/p_{_{2}}}\\
&&\; \le C \sup_{k_{0}\in \mathbb{Z}}2^{-k_{0}\lambda p_{_{1}}}  \bigg(\sum_{k=-\infty}^{k_{0}}2^{k\alpha p_{_{1}}}
  \|I_{\beta, b}^{m}(f)\cdot\chi_{_{\scriptstyle{k}}}  \|_{L^{q_{_{2}}(\cdot)}(\R^{n})}^{p_{_{1}}}\bigg) \\
&&\; \le C \sup_{k_{0}\in \mathbb{Z}}2^{-k_{0}\lambda p_{_{1}}}  \bigg(\sum_{k=-\infty}^{k_{0}}2^{k\alpha p_{_{1}}}
  \Big(\sum_{j=-\infty}^{k-2}\|I_{\beta, b}^{m}(f_{j})\cdot\chi_{_{\scriptstyle{k}}}  \|_{L^{q_{_{2}}(\cdot)}(\R^{n})} \Big)^{p_{_{1}}}\bigg) \\
&&\;\ \ + C \sup_{k_{0}\in \mathbb{Z}}2^{-k_{0}\lambda p_{_{1}}}  \bigg(\sum_{k=-\infty}^{k_{0}}2^{k\alpha p_{_{1}}}
  \Big(\sum_{j=k-1}^{k+1}\|I_{\beta, b}^{m}(f_{j})\cdot\chi_{_{\scriptstyle{k}}}  \|_{L^{q_{_{2}}(\cdot)}(\R^{n})} \Big)^{p_{_{1}}}\bigg)\\
&&\;\ \ + C \sup_{k_{0}\in \mathbb{Z}}2^{-k_{0}\lambda p_{_{1}}}  \bigg(\sum_{k=-\infty}^{k_{0}}2^{k\alpha p_{_{1}}}
  \Big(\sum_{j=k+2}^{\infty}\|I_{\beta, b}^{m}(f_{j})\cdot\chi_{_{\scriptstyle{k}}}  \|_{L^{q_{_{2}}(\cdot)}(\R^{n})} \Big)^{p_{_{1}}}\bigg)\\
 &&\;  = C(E_{1} +E_{2}+E_{3}).
\end{eqnarray*}

First we estimate $E_{2}$. Using the Proposition \ref{Pr-2}, we have
\begin{eqnarray*}
&&\; E_{2} =  \sup_{k_{0}\in \mathbb{Z}}2^{-k_{0}\lambda p_{_{1}}}  \bigg(\sum_{k=-\infty}^{k_{0}}2^{k\alpha p_{_{1}}}
  \Big(\sum_{j=k-1}^{k+1}\|I_{\beta, b}^{m}(f_{j})\cdot\chi_{_{\scriptstyle{k}}}  \|_{L^{q_{_{2}}(\cdot)}(\R^{n})} \Big)^{p_{_{1}}}\bigg)\\
&&\;\le C\|b\|^{mp_{1}}_{\BMO(\R^{n})} \sup_{k_{0}\in \mathbb{Z}}2^{-k_{0}\lambda p_{_{1}}}  \bigg(\sum_{k=-\infty}^{k_{0}}2^{k\alpha p_{_{1}}}
  \Big(\sum_{j=k-1}^{k+1}\|f_{j}\cdot\chi_{_{\scriptstyle{k}}}  \|_{L^{q_{_{1}}(\cdot)}(\R^{n})} \Big)^{p_{_{1}}}\bigg)\\
  &&\;\le C\|b\|^{mp_{1}}_{\BMO(\R^{n})} \sup_{k_{0}\in \mathbb{Z}}2^{-k_{0}\lambda p_{_{1}}}  \bigg(\sum_{k=-\infty}^{k_{0}}2^{k\alpha p_{_{1}}}\|f_{j}\cdot\chi_{_{\scriptstyle{k}}}  \|_{L^{q_{_{1}}(\cdot)}(\R^{n})} ^{p_{_{1}}}\bigg)\\
 &&\;  = C\|b\|^{mp_{1}}_{\BMO(\R^{n})} \|f\|^{p_{1}}_{M\dot{K}_{p_{_{1}},q_{_{1}}(\cdot)}^{\alpha,\lambda}(\R^{n})}.
\end{eqnarray*}

For $E_{1}$. Note that when $x\in A_{k}, j\le k-2$, and $y\in A_{j}$, then $|x-y|\backsim |x|, 2|y|\le |x|$. Therefore, using the generalized H\"{o}lder's inequality(see 1), Lemma \ref{Lem.1}), we have
\begin{equation*}
\begin{split}
& |I_{\beta, b}^{m}(f_{j})(x)\cdot\chi_{_{\scriptstyle k}}(x)|  \le C \int_{A_{j}}\frac{|f_{j}(y)| |b(x)-b(y)|^{m}}{|x-y|^{n-\beta}} \mathrm{d}y \cdot\chi_{_{\scriptstyle k}}(x)   \\
& \le C2^{^{k(\beta-n)}} \int_{A_{j}}|f_{j}(y)||b(x)-b(y)|^{m}   \mathrm{d}y \cdot\chi_{_{\scriptstyle{k}}}(x)\\
& \le C2^{^{k(\beta-n)}} \bigg(|b(x)-b_{B_{j}}|^{m} \int_{A_{j}}|f_{j}(y)| \mathrm{d}y+\int_{A_{j}}|f_{j}(y)||b(y)-b_{B_{j}}|^{m}   \mathrm{d}y \bigg) \cdot\chi_{_{\scriptstyle{k}}}(x)\\
& \le C2^{^{k(\beta-n)}} \|f_{j}\|_{L^{^{q_{_{1}}(\cdot)}}(\R^{n})} \bigg(|b(x)-b_{B_{j}}|^{m} \|\chi_{_{\scriptstyle{j}}}\|_{L^{^{q'_{_{1}}(\cdot)}}(\R^{n})} \\
 &\ \ \hspace{10cc} +\|(b-b_{B_{j}})^{m} \chi_{_{\scriptstyle{j}}}\|_{L^{^{q'_{_{1}}(\cdot)}}(\R^{n})} \bigg) \cdot\chi_{_{\scriptstyle{k}}}(x).
\end{split}
\end{equation*}

Thus, from Lemma \ref{Lem.4}, and note that $\|\chi_{_{\scriptstyle{i}}}\|_{L^{^{s(\cdot)}}(\R^{n})} \le\|\chi_{_{\scriptstyle{B_{i}}}}\|_{L^{^{s(\cdot)}}(\R^{n})}$, it follows that
\begin{equation}
\begin{split}
& \|I_{\beta, b}^{m}(f_{j})\cdot\chi_{_{\scriptstyle k}}\|_{L^{q_{_{2}}(\cdot)}(\R^{n})}  \le C2^{^{k(\beta-n)}} \|f_{j}\|_{L^{^{q_{_{1}}(\cdot)}}(\R^{n})} \bigg(\|(b-b_{B_{j}})^{m}\chi_{_{\scriptstyle{k}}}\|_{L^{q_{_{2}}(\cdot)}(\R^{n})}\|\chi_{_{\scriptstyle{j}}}\|_{L^{^{q'_{_{1}}(\cdot)}}(\R^{n})} \\ &\ \ \hspace{7cc} +\|(b-b_{B_{j}})^{m}\chi_{_{\scriptstyle{j}}}\|_{L^{^{q'_{_{1}}(\cdot)}}(\R^{n})}\|\chi_{_{\scriptstyle{k}}}\|_{L^{q_{_{2}}(\cdot)}(\R^{n})}\bigg)   \\
& \le C2^{^{k(\beta-n)}} \|f_{j}\|_{L^{^{q_{_{1}}(\cdot)}}(\R^{n})} \bigg((k-j)^{m} \|b\|^{m}_{\BMO(\R^{n})} \| \chi_{B_{k}}\|_{L^{^{q_{_{2}}(\cdot)}}(\R^{n})} \|\chi_{_{\scriptstyle{j}}}\|_{L^{^{q'_{_{1}}(\cdot)}}(\R^{n})} \\
&\ \ \hspace{7cc}  +\|b\|^{m}_{\BMO(\R^{n})} \|\chi_{B_{j }}\|_{L^{^{q'_{_{1}}(\cdot)}}(\R^{n})} \|\chi_{_{\scriptstyle{k}}}\|_{L^{q_{_{2}}(\cdot)}(\R^{n})}\bigg)   \\
& \le C2^{^{k(\beta-n)}}  (k-j)^{m}\|b\|^{m}_{\BMO(\R^{n})}  \|f_{j}\|_{L^{^{q_{_{1}}(\cdot)}}(\R^{n})} \|\chi_{B_{j}}\|_{L^{^{q'_{_{1}}(\cdot)}}(\R^{n})}  \| \chi_{B_{k}}\|_{L^{^{q_{_{2}}(\cdot)}}(\R^{n})}.
\end{split}      \label{H-point}
\end{equation}

Note that $\chi_{_{\scriptstyle{B_{k}}}}(x)\le {C}2^{-k\beta}
I_{\beta}(\chi_{_{\scriptstyle{B_{k}}}})(x)$~(see page 350,\cite{I}),  by Proposition \ref{Pr-3} and Lemma \ref{Lem.3}, we obtain
\begin{eqnarray}
\begin{split} \label{F-K}
\|\chi_{_{\scriptstyle B_{k}}}\|_{L^{^{q_{_{2}}(\cdot)}}(\R^{n})}     & \le C 2^{^{-k\beta}} \|I_{\beta}(\chi_{_{\scriptstyle{B_{k}}}}) \|_{L^{^{q_{_{2}}(\cdot)}}(\R^{n})} \\
& \le C2^{^{-k\beta}} \|\chi_{_{\scriptstyle B_{k}}} \|_{L^{^{q_{_{1}}(\cdot)}}(\R^{n})}.
 \end{split}
\end{eqnarray}
Using Lemma \ref{Lem.2}, Lemma \ref{Lem.3}, (\ref{BZ})  and (\ref{F-K}), we have
\begin{eqnarray} \label{H1-norm}
&&\; 2^{^{k(\beta-n)}} \|\chi_{B_{j}}\|_{L^{^{q'_{_{1}}(\cdot)}}(\R^{n})} \|\chi_{B_{k}}\|_{L^{^{q_{_{2}}(\cdot)}}(\R^{n})} \le 2^{^{k(\beta-n)}} \|\chi_{B_{j}}\|_{L^{^{q'_{_{1}}(\cdot)}}(\R^{n})} \cdot 2^{^{-k\beta}} \|\chi_{_{\scriptstyle B_{k}}} \|_{L^{^{q_{_{1}}(\cdot)}}(\R^{n})} \nonumber\\
&&\; \le C  \|\chi_{B_{j}}\|_{L^{^{q'_{_{1}}(\cdot)}}(\R^{n})}  \cdot 2^{^{-kn}} \|\chi_{_{\scriptstyle B_{k}}} \|_{L^{^{q_{_{1}}(\cdot)}}(\R^{n})}     \le C  \|\chi_{B_{j}}\|_{L^{^{q'_{_{1}}(\cdot)}}(\R^{n})} \|\chi_{_{\scriptstyle B_{k}}} \|^{-1}_{L^{^{q'_{_{1}}(\cdot)}}(\R^{n})}     \\
&&\;  = C\frac{\|\chi_{_{\scriptstyle {B_{j}}}} \|_{L^{^{q'_{_{1}}(\cdot)}}(\R^{n})}}
 {\|\chi_{_{\scriptstyle{B_{k}}}}\|_{L^{^{q'_{_{1}}(\cdot)}}(\R^{n})}} \le C 2^{(j-k)n\delta_{1}}. \nonumber
  \end{eqnarray}

On the other hand, note the following fact
\begin{eqnarray}
\begin{split} \label{H1-f-j}
\|f_{j}\|_{L^{^{q_{_{1}}(\cdot)}}(\R^{n})} &= 2^{-j\alpha}\Big(2^{j{\alpha}p_{_1}} \|f_{j}\|^{p_{_1}}_{L^{^{q_{_{1}}(\cdot)}}(\R^{n})}\Big)^{1/p_{_1}}\\
&\le  2^{-j\alpha}\bigg(\sum_{i=-\infty}^j 2^{i{\alpha}p_{_1}} \|f_{i}\|^{p_{_1}}_{L^{^{q_{_{1}}(\cdot)}}(\R^{n})}\bigg)^{1/p_{_1}}\\
&=  2^{j(\lambda-\alpha)}\bigg(2^{-j\lambda} \Big(\sum_{i=-\infty}^j 2^{i{\alpha}p_{_1}}
 \|f_{i}\|^{p_{_1}}_{L^{^{q_{_{1}}(\cdot)}} (\R^{n})}\Big)^{1/p_{_1}}\bigg)\\
&\le  C 2^{j(\lambda-\alpha)} \|f\|_{M\dot{K}_{p_{_{1}},q_{_{1}}(\cdot)}^{\alpha,\lambda}(\R^{n})}.
 \end{split}
\end{eqnarray}

Thus, combing (\ref{H-point}), (\ref{H1-norm})  and (\ref{H1-f-j}), and using $\alpha< \lambda+n\delta_{1}$, it follows that
\begin{equation*}
\begin{split}
E_{1} &= \sup_{k_{0}\in \mathbb{Z}}2^{-k_{0}\lambda p_{_{1}}}  \bigg(\sum_{k=-\infty}^{k_{0}}2^{k\alpha p_{_{1}}}
  \Big(\sum_{j=-\infty}^{k-2}\|I_{\beta, b}^{m}(f_{j})\cdot\chi_{_{\scriptstyle{k}}}  \|_{L^{q_{_{2}}(\cdot)}(\R^{n})} \Big)^{p_{_{1}}}\bigg) \\
&\le C\sup_{k_{0}\in \mathbb{Z}}2^{-k_{0}\lambda p_{_{1}}}  \bigg(\sum_{k=-\infty}^{k_{0}}2^{k\alpha p_{_{1}}}
  \Big(\sum_{j=-\infty}^{k-2}(k-j)^{m}\|b\|^{m}_{\BMO(\R^{n})}  \|f_{j}\|_{L^{^{q_{_{1}}(\cdot)}}(\R^{n})}  2^{-(k-j)n\delta_{1}} \Big)^{p_{_{1}}}\bigg) \\
&\le C\|b\|^{mp_{1}}_{\BMO(\R^{n})} \|f\|^{p_{1}}_{M\dot{K}_{p_{_{1}},q_{_{1}}(\cdot)}^{\alpha,\lambda}(\R^{n})} \\
&\ \hspace{2cc} \times \sup_{k_{0}\in \mathbb{Z}}2^{-k_{0}\lambda p_{_{1}}}  \bigg(\sum_{k=-\infty}^{k_{0}}2^{k\lambda p_{_{1}}}  \Big(\sum_{j=-\infty}^{k-2}(k-j)^{m} 2^{(k-j)(\alpha-\lambda-n\delta_{1})} \Big)^{p_{_{1}}}\bigg) \\
&\le C\|b\|^{mp_{1}}_{\BMO(\R^{n})} \|f\|^{p_{1}}_{M\dot{K}_{p_{_{1}},q_{_{1}}(\cdot)}^{\alpha,\lambda}(\R^{n})} \sup_{k_{0}\in \mathbb{Z}}2^{-k_{0}\lambda p_{_{1}}}  \bigg(\sum_{k=-\infty}^{k_{0}}2^{k\lambda p_{_{1}}}\bigg) \\
&\le C\|b\|^{mp_{1}}_{\BMO(\R^{n})} \|f\|^{p_{1}}_{M\dot{K}_{p_{_{1}},q_{_{1}}(\cdot)}^{\alpha,\lambda}(\R^{n})}.
\end{split}
\end{equation*}

Now, let us turn to estimate for $E_{3}$. Note that when $x\in A_{k}, j\ge k+2$, and $y\in A_{j}$, then $|x-y|\backsim |y|, 2|x|\le |y|$. Therefore, using the generalized H\"{o}lder's inequality(see 1), Lemma \ref{Lem.1}), we have
\begin{equation*}
\begin{split}
& |I_{\beta, b}^{m}(f_{j})(x)\cdot\chi_{_{\scriptstyle k}}(x)|  \le C \int_{A_{j}}\frac{|f_{j}(y)| |b(x)-b(y)|^{m}}{|x-y|^{n-\beta}} \mathrm{d}y \cdot\chi_{_{\scriptstyle k}}(x)   \\
& \le C2^{^{j(\beta-n)}} \int_{A_{j}}|f_{j}(y)||b(x)-b(y)|^{m}   \mathrm{d}y \cdot\chi_{_{\scriptstyle{k}}}(x)\\
& \le C2^{^{j(\beta-n)}} \bigg(|b(x)-b_{B_{k}}|^{m} \int_{A_{j}}|f_{j}(y)| \mathrm{d}y+\int_{A_{j}}|f_{j}(y)||b(y)-b_{B_{k}}|^{m}   \mathrm{d}y \bigg) \cdot\chi_{_{\scriptstyle{k}}}(x)\\
& \le C2^{^{j(\beta-n)}} \|f_{j}\|_{L^{^{q_{_{1}}(\cdot)}}(\R^{n})} \bigg(|b(x)-b_{B_{k}}|^{m} \|\chi_{_{\scriptstyle{j}}}\|_{L^{^{q'_{_{1}}(\cdot)}}(\R^{n})} +\|(b-b_{B_{k}})^{m} \chi_{_{\scriptstyle{j}}}\|_{L^{^{q'_{_{1}}(\cdot)}}(\R^{n})} \bigg) \cdot\chi_{_{\scriptstyle{k}}}(x).
\end{split}      
\end{equation*}

Using Lemma \ref{Lem.4}, it follows that
\begin{eqnarray} \label{H-point3}
&&\; \|I_{\beta, b}^{m}(f_{j})\cdot\chi_{_{\scriptstyle k}}\|_{L^{q_{_{2}}(\cdot)}(\R^{n})}
 \le C2^{^{j(\beta-n)}} \|f_{j}\|_{L^{^{q_{_{1}}(\cdot)}}(\R^{n})} \bigg(\|(b-b_{B_{k}})^{m}\chi_{_{\scriptstyle{k}}}\|_{L^{q_{_{2}}(\cdot)}(\R^{n})}\|\chi_{_{\scriptstyle{j}}}\|_{L^{^{q'_{_{1}}(\cdot)}}(\R^{n})} \nonumber\\
&&\;\ \ \hspace{10cc} +\|(b-b_{B_{k}})^{m}\chi_{_{\scriptstyle{j}}}\|_{L^{^{q'_{_{1}}(\cdot)}}(\R^{n})}\|\chi_{_{\scriptstyle{k}}}\|_{L^{q_{_{2}}(\cdot)}(\R^{n})}\bigg)   \nonumber\\
&&\; \le C2^{^{j(\beta-n)}} \|f_{j}\|_{L^{^{q_{_{1}}(\cdot)}}(\R^{n})} \bigg( \|b\|^{m}_{\BMO(\R^{n})} \| \chi_{B_{k}}\|_{L^{^{q_{_{2}}(\cdot)}}(\R^{n})} \|\chi_{_{\scriptstyle{j}}}\|_{L^{^{q'_{_{1}}(\cdot)}}(\R^{n})} \\
&&\;\ \ \hspace{10cc}  +(j-k)^{m}\|b\|^{m}_{\BMO(\R^{n})} \|\chi_{B_{j }}\|_{L^{^{q'_{_{1}}(\cdot)}}(\R^{n})} \|\chi_{_{\scriptstyle{k}}}\|_{L^{q_{_{2}}(\cdot)}(\R^{n})}\bigg)    \nonumber\\
&&\; \le C2^{^{j(\beta-n)}}  (j-k)^{m}\|b\|^{m}_{\BMO(\R^{n})}  \|f_{j}\|_{L^{^{q_{_{1}}(\cdot)}}(\R^{n})} \|\chi_{B_{j}}\|_{L^{^{q'_{_{1}}(\cdot)}}(\R^{n})}  \| \chi_{B_{k}}\|_{L^{^{q_{_{2}}(\cdot)}}(\R^{n})}.   \nonumber
\end{eqnarray}

Note that $\chi_{_{\scriptstyle{B_{j}}}}(x)\le {C}2^{-j\beta}
I_{\beta}(\chi_{_{\scriptstyle{B_{j}}}})(x)$~(see page 350,\cite{I}),  by Proposition \ref{Pr-3} and Lemma \ref{Lem.3}, we obtain
\begin{eqnarray*}
\begin{split}
\|\chi_{B_{j}}\|_{L^{^{q_{_{2}}(\cdot)}}(\R^{n})} & \le C2^{^{-j\beta}} \|I_{\beta}(\chi_{_{\scriptstyle{B_{j}}}}) \|_{L^{^{q_{_{2}}(\cdot)}}(\R^{n})} \\
& \le C2^{^{-j\beta}} \|\chi_{_{\scriptstyle B_{j}}} \|_{L^{^{q_{_{1}}(\cdot)}}(\R^{n})}   \\
& \le C2^{^{-j\beta}} 2^{^{jn}}\|\chi_{_{\scriptstyle B_{j}}} \|^{-1}_{L^{^{q'_{_{1}}(\cdot)}}(\R^{n})}.
 \end{split}
\end{eqnarray*}
Thus, we have
\begin{eqnarray}
\begin{split} \label{F-J}
2^{^{j(\beta-n)}}\|\chi_{_{\scriptstyle B_{j}}} \|_{L^{^{q'_{_{1}}(\cdot)}}(\R^{n})} \le C \|\chi_{B_{j}}\|^{-1}_{L^{^{q_{_{2}}(\cdot)}}(\R^{n})}.
 \end{split}
\end{eqnarray}

Using Lemma \ref{Lem.2}, Lemma \ref{Lem.3}, (\ref{BZ})  and (\ref{F-J}), we have
\begin{eqnarray}
\begin{split} \label{H1-norm3}
2^{^{j(\beta-n)}} \|\chi_{B_{j}}\|_{L^{^{q'_{_{1}}(\cdot)}}(\R^{n})}\|\chi_{B_{k}}\|_{L^{^{q_{_{2}}(\cdot)}}(\R^{n})} &\le C \|\chi_{B_{j}}\|^{-1}_{L^{^{q_{_{2}}(\cdot)}}(\R^{n})} \|\chi_{B_{k}}\|_{L^{^{q_{_{2}}(\cdot)}}(\R^{n})} \\
 &\le C \frac{\|\chi_{_{\scriptstyle {B_{k}}}} \|_{L^{^{q_{_{2}}(\cdot)}}(\R^{n})}}{\|\chi_{B_{j}}\|_{L^{^{q_{_{2}}(\cdot)}}(\R^{n})} }\\
     & \le C  2^{(k-j)n\delta_{2}}.
 \end{split}
\end{eqnarray}

Thus, combing (\ref{H1-f-j}),(\ref{H-point3}) and(\ref{H1-norm3}), and using $\lambda-n\delta_{2}<\alpha$, it follows that
\begin{eqnarray*}
\begin{split}
E_{3} &= \sup_{k_{0}\in \mathbb{Z}}2^{-k_{0}\lambda p_{_{1}}}  \bigg(\sum_{k=-\infty}^{k_{0}}2^{k\alpha p_{_{1}}}
  \Big(\sum_{j=k+2}^{\infty}\|I_{\beta, b}^{m}(f_{j})\cdot\chi_{_{\scriptstyle{k}}}  \|_{L^{q_{_{2}}(\cdot)}(\R^{n})} \Big)^{p_{_{1}}}\bigg) \\
&\le C\sup_{k_{0}\in \mathbb{Z}}2^{-k_{0}\lambda p_{_{1}}}  \bigg(\sum_{k=-\infty}^{k_{0}}2^{k\alpha p_{_{1}}}
  \Big(\sum_{j=k+2}^{\infty}(j-k)^{m}\|b\|^{m}_{\BMO(\R^{n})}  \|f_{j}\|_{L^{^{q_{_{1}}(\cdot)}}(\R^{n})}  2^{-(j-k)n\delta_{2}} \Big)^{p_{_{1}}}\bigg) \\
&\le C\|b\|^{mp_{1}}_{\BMO(\R^{n})} \|f\|^{p_{1}}_{M\dot{K}_{p_{_{1}},q_{_{1}}(\cdot)}^{\alpha,\lambda}(\R^{n})} \\
&\ \hspace{2cc} \times \sup_{k_{0}\in \mathbb{Z}}2^{-k_{0}\lambda p_{_{1}}}  \bigg(\sum_{k=-\infty}^{k_{0}}2^{k\lambda p_{_{1}}}  \Big(\sum_{j=k+2}^{\infty}(j-k)^{m}2^{(j-k)(\lambda-\alpha-n\delta_{2})} \Big)^{p_{_{1}}}\bigg) \\
&\le C\|b\|^{mp_{1}}_{\BMO(\R^{n})} \|f\|^{p_{1}}_{M\dot{K}_{p_{_{1}},q_{_{1}}(\cdot)}^{\alpha,\lambda}(\R^{n})} \sup_{k_{0}\in \mathbb{Z}}2^{-k_{0}\lambda p_{_{1}}}  \bigg(\sum_{k=-\infty}^{k_{0}}2^{k\lambda p_{_{1}}}\bigg) \\
&\le C\|b\|^{mp_{1}}_{\BMO(\R^{n})} \|f\|^{p_{1}}_{M\dot{K}_{p_{_{1}},q_{_{1}}(\cdot)}^{\alpha,\lambda}(\R^{n})}.
\end{split}
\end{eqnarray*}

This finishes the proof of Theorem \ref{thm.1}.
\end{proof}

When $\lambda=0$, our main result also hold on Herz space with variable expoent, and generalize the result of Izuki\upcite{I3} (see Theorem 3).  When $m=0$, we also improve the result for Izuki\upcite{I} (see Theorem 2).

\begin{acknowledgments}
 The author cordially  thank the referees for their valuable suggestions and useful comments which have lead to the improvement of this paper. This work was supported  in part  by  the Pre-Research Project (No. SY201224) of Provincial Key Innovation  and NSF (No. 11161042)  of China.
\end{acknowledgments}

\end{document}